\documentclass[10pt]{article}

\usepackage[all]{xy}
\usepackage{graphics}
\usepackage{color}
\usepackage{url}
\usepackage{bm}
\usepackage[tbtags]{amsmath}
\usepackage{latexsym,array,delarray,amsfonts,amssymb,amsthm}

\theoremstyle{plain}
\newtheorem{thm}{Theorem}[section]
\newtheorem{lmma}[thm]{Lemma}

\newtheorem{cor}[thm]{Corollary}

\newtheorem{fact}[thm]{Fact}

\newtheorem{ex}[thm]{Example}

\newtheorem{rmk}[thm]{Remark}

\newcommand{\Q}{{\mathbb Q}}
\newcommand{\Z}{{\mathbb Z}}
\newcommand{\N}{{\mathbb N}}

\newcommand{\C}{{\mathbb C}}
\newcommand{\R}{\mathbb R}

\begin{document}

\medskip
\title{\bf On the Enumeration of Certain Weighted Graphs}

\author{ Mikl\'os B\'ona\thanks{ {\tt bona@math.ufl.edu},
Department of Mathematics, University of Florida, Gainesville, FL
32611 USA.} \and Hyeong-Kwan Ju\thanks{ {\tt hkju@chonnam.ac.kr},
Department of Mathematics, Chonnam National University, Kwangju,
500-757, Republic of Korea. The research of this author was
financially supported by Chonnam National University. } \and
Ruriko Yoshida\thanks{{\tt ruriko@math.duke.edu}, Department of
Mathematics, Duke University, Durham, NC 27708-0320 USA. }}

\date{\today}

\maketitle
\begin{abstract}  We enumerate  weighted graphs with a certain
upper bound condition. We also compute the generating function of
the numbers of these graphs, and prove that it is a rational
function. In particular, we show that if the given graph is a
bipartite graph, then its generating function is of the form
$\frac{p(x)}{(1-x)^{m+1}}$, where $m$ is the number of vertices of
the graph and $p(x)$ is a polynomial of degree at most $m$.
\end{abstract}


\section{Introduction}\label{intro}

For a given nonnegative integer $n$, let $[n]:=\{0,1,2,\cdots ,
n\}$. Also, let $G=(V,E)$ be a simple graph (no loops and no
multiple edges allowed) with the vertex set $V= \{v_1, v_2, \cdots
, v_m \}$. Let $ \alpha =( n_1, n_2, \cdots , n_m ) \in [n]^m$ so
that

\begin{equation}\label{eq:1}
 n_i + n_j \le n \quad \mbox{ if } \quad
v_i v_j \in E. \end{equation} In other words, the sum of two
weights corresponding to an adjacent pair of vertices is bounded
by a given integer $n$.

We call a triplet $WG_\alpha =(V,E,\alpha )$ a {\bf weighted
graph} of G with  distribution $\alpha$. We will denote by $WG(n)$
the number of  all weighted graphs of $G$ with a fixed upper bound
$n$. Note that $WG(n)$ enumerates the set of all solutions of the
system of linear inequalities corresponding to the graph $G$ by
the condition \eqref{eq:1}.\\

Let ${\mathcal G}$ be the set of all simple graphs and ${\bf
C}[[x]]$ the ring of formal power series. Define a map $$ \rho :
{\mathcal G} \rightarrow \C[[x]]\ (G\ \mapsto \rho(G)=
\sum_{n=0}^\infty WG(n)x^n).$$

For example, consider the following graph $G$: \xymatrix{ a
\ar@{-}[r]
\ar@{-}[dr] \ar@{-}[d] & b \ar@{-}[r]  & c \\
d & e }.

For each $a=i \in [n],\hspace{0.3cm} b,d,e $ can be $0,1,2,\cdots,
n-i$, and for each $b=j \in [n-i]$, $c$  can be
$0,1,2,\cdots,n-j$. Hence,
\begin{equation*}\begin{aligned}
WG(n) &= \sum_{i=0}^{n}(n+1-i)^{2}(\sum_{j=0}^{n-i}(n+1-j))\\
        &= \frac{(3n+5)(2n+3)(3n+4)(n+2)(n+1)}{120}\\
        &= \frac{1}{120} (18(n+1)^5 +45(n+1)^4 + 40(n+1)^3  + 15(n+1)^2
+2(n+1)),
\end{aligned}\end{equation*}
and
\begin{equation*}\begin{aligned}
\rho(G) &= \frac{1}{120} (18 \frac{1+26x+66x^2 +26x^3 +
x^4}{(1-x)^6} +
45 \frac{1+11x+11x^2 +x^3}{(1-x)^5}\\
        &+ 40\frac{1+4x+x^2}{(1-x)^4} + 15 \frac{1+x}{(1-x)^3} + 2
        \frac{1}{(1-x)^2})\\
        &= \frac{(1+x)(1+7x+x^2)}{(1-x)^6} =\frac{1+8x+8x^2+x^3}{(1-x)^6}.
\end{aligned}\end{equation*}

Equivalently, $WG(n)$ is the number of solutions $(n_a, n_b, n_c,
n_d, n_e) \in [n]^5 $ to the following system of inequalities:
$$ \left\{ \begin{array}{c}
n_a + n_b \le n\\
n_b + n_c \le n\\
n_a + n_d \le n\\
n_a + n_e \le n
\end{array} \right.
$$

\begin{rmk}

\bigskip

\noindent
\begin{enumerate}
\item The map $\rho$ is not injective.\\
 For graphs $G_1$ 
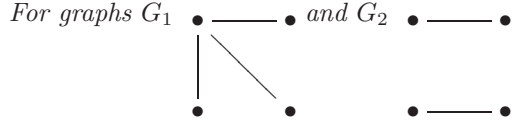
\xymatrix{
{\bullet} \ar@{-}[r]
\ar@{-}[dr] \ar@{-}[d] & {\bullet} \\
{\bullet} & {\bullet} } and  $G_2$ \xymatrix{ {\bullet} \ar@{-}[r]
& {\bullet} \\
 {\bullet} \ar@{-}[r] & {\bullet} } \quad
 $$\rho(G_1) = \rho(G_2)=
 \frac{1+4x+x^2}{(1-x)^5} \quad
\mbox{since} \quad \sum_{k=1}^{n+1} k^3 = (\sum_{k=1}^{n+1}
k)^2.$$
\item We used the following identities (to get $\rho(G)$ in the
previous example):
$$(n+1)^m = \sum_{k=1}^n A(m,k) C(n+k,m)$$ where $A(m,k)$ is
an Eulerian number, and $C(n,m)=\frac{n!}{m!(n-m)!}$ is a binomial
coefficient (See Fact \ref{Facts} 3. in Section \ref{disGraph}.),
and

 $$ \sum_{n=0}^{\infty} C(n+k,k) x^n = \frac{1}{(1-x)^{k+1}}.
$$
\end{enumerate}
\end{rmk}

 In Section \ref{genGraph} we compute $\rho(G)$ for linear graphs,
circular graphs, complete graphs, star graphs, discrete graphs, a
cubic graph, an octahedral graph and complete bipartite graphs. In
Section \ref{genFun} we remind a reader of rational convex
polytopes and rational generating functions, and describe the
relationship between these objects and our problem (enumeration of
weighted graphs). In Section \ref{question} we list many new
problems that we encountered while we worked on this topic.

\section{Generating Functions for Various Graphs}\label{genGraph}
\subsection{Linear Graphs}\label{linGraph}

For the null graph $\emptyset$ we set $W\emptyset(n)=1$ for all
nonnegative integers $n$, so $$\rho(\emptyset)= \frac{1}{1-x}.$$
Let $L_1$ be a one-vertex graph, and for $i \geq 2$ let $L_i$ be
the linear graph with $i$ vertices. That is, $L_i$ is a tree with
two vertices of degree one and $i-2$ vertices of degree of two.
For the one-vertex graph $L_1=K_1 (\bullet), \quad WL_1(n)=n+1$
and
$$\rho(L_1)= \frac{1}{(1-x)^2}.$$

For $L_2=K_2$ (\xymatrix{a{\bullet} \ar@{-}[r] & {\bullet}b\\
&}) let $a=i,\ b=j$, and

$$b_{ij}= \begin{cases}
         1 &\text{if } i+j \le n,\\
         0 &\text{otherwise}
        \end{cases}
       \quad = \quad \lfloor{\frac{2n+3-i-j}{n+2}}\rfloor .$$
Note that all indices $i$ and $j$ run from $0$ to $n$, not from
$1$ to $n+1$. $B(n):=(b_{ij})$ is an $(n+1)\times(n+1)$ matrix,
$J_n: =(1,1,1,\cdots, 1)^t \in [n]^{n+1}$. Then
$$WL_2(n)=J_n^t B(n) J_n= \sum_{i,j=0}^n b_{ij}=C(n+2,2),$$
$$\rho(L_2)= \frac{1}{(1-x)^3}.$$

For $L_3$ (\xymatrix{{\bullet} \ar@{-}[r] & {\bullet} \ar@{-}[r] & {\bullet}\\
&&}),
$$WL_3(n)=J_n^t B(n)^2 J_n=C(n+3,3)+C(n+2,3),$$
$$\rho(L_3)=\frac{1+x}{(1-x)^4}.$$
By induction on $k$ we can show that
$$WL_{k+1}(n)=J_n^t B(n)^k J_n$$ holds for all nonnegative integers $k$.

For example,\\

$B(4) =\begin{bmatrix}
1 & 1 & 1 & 1 & 1 \\
1 & 1 & 1 & 1 & 0 \\
1 & 1 & 1 & 0 & 0 \\
1 & 1 & 0 & 0 & 0 \\
1 & 0 & 0 & 0 & 0
\end{bmatrix},\quad
{B(4)}^2 = \begin{bmatrix}
5 & 4 & 3 & 2 & 1 \\
4 & 4 & 3 & 2 & 1 \\
3 & 3 & 3 & 2 & 1 \\
2 & 2 & 2 & 2 & 1 \\
1 & 1 & 1 & 1 & 1
\end{bmatrix},$
$$$$
${B(4)}^3 = \begin{bmatrix}
15 & 14 & 12 & 9 & 5 \\
14 & 13 & 11 & 8 & 4 \\
12 & 11 &  9 & 6 & 3 \\
9 &  8 &  6 & 4 & 2 \\
5 &  4 &  3 & 2 & 1
\end{bmatrix},\quad
{B(4)}^4 = \begin{bmatrix}
55 & 50 & 41 & 29 & 15 \\
50 & 46 & 38 & 27 & 14 \\
41 & 38 & 32 & 23 & 12 \\
29 & 27 & 23 & 17 & 9 \\
15 & 14 & 12 & 9 & 5
\end{bmatrix}.$

\medskip
Note that $WL_{k+1}(n)$ is the $(0,0)$-entry of $B(n)^{k+2}$. For
example, $WL_3(4)=55,$ which is the sum of all entries of the
matrix $B(4)^2.$
\bigskip

\begin{thm}[Bona and Ju, 2006 \cite{BJ}]
 Let $$F(n,x)=\sum_{k=0}^\infty WL_{k+1}(n)x^k =
\sum_{k=0}^\infty (J_n^t B(n)^k J_n)x^k.$$ Then
$$F(n,x)=\frac{1+F(n-1,-x)}{1-x(F(n-1,-x)+1)}
      =\cfrac{1}{-x + \cfrac{1}{1+F(n-1,-x)}}$$
     $$ = \cfrac{1}{-x + \cfrac{1}{1 + \cfrac{1}{x + \cfrac{1}{1 +
F(n-2,x)}}}}
      :=[-x,1,x,1,F(n-2,x)],$$
where $F(0,x)=\frac{1}{-x+1}=[-x,1]$ and
$F(1,0)=\frac{2+x}{1-x-x^2}=[-x,1,x,1].  $
\end{thm}
\bigskip

\begin{thm}[Bona and Ju, 2006 \cite{BJ}]
 $$F(n,x)=\frac{P_n(x)}{Q_n(x)},$$
where $$Q_n(x)=\det(I-xB(n))= \sum_{k=0}^n \binom{\lfloor
\frac{n+k+1}{2}\rfloor}{k} (-1)^{\lfloor \frac{k+1}{2} \rfloor}
x^k$$ and $$P_n(x)=(Q_{n-2}(x)-(1+x)Q_n(x))/x^2$$ is a polynomial
of degree $deg(Q_n(x))-1=n.\quad  $
\end{thm}

\bigskip

\begin{rmk}

\bigskip

\noindent
\begin{enumerate}
\item $Q_n(x)$ satisfies a recurrence relation
\begin{equation}\label{eq:222}
Q_n(x)+(x^2-2)Q_{n-2}(x)+Q_{n-4}(x)=0.
\end{equation}

\item From the recurrence relation \eqref{eq:222} we get the following
generating function $R(x,y)$ for the sequence
$\{Q_n(x)\}_{n=0}^{\infty}$:
$$R(x,y)=\sum_{n=0}^\infty Q_{n-1}(x)y^n =
\frac{(1+y)(1-xy-y^2)}{1-y^2(2-x^2-y^2)}.$$

\item $\{F(n,x)\}_{n=0}^\infty$ is an approximation(called P\'ade
approximant) to the infinite periodic (of period 4) continued
fraction $$[-x,1,x,1,-x,1,x,1, \cdots]\\
=\frac{-(2+x) \pm \sqrt{x^2-4}}{2x},$$ and converges to
$\begin{cases}\frac{-(2+x) + \sqrt{x^2-4}}{2x}  &{\mbox on }\quad
[2, \infty),\\ \frac{-(2+x) - \sqrt{x^2-4}}{2x}  &{\mbox on }\quad
(-\infty, 2]. \end{cases}$

\item The explicit expression of the generating function $G(x,y)$ for
the sequence $\{F(n,x)\}_{n=0}^{\infty}$ is not known, but from
Theorem 2.2
$$G(x,y)=\sum_{n=0}^\infty F(n,x)y^n = \sum_{n=0}^\infty
\frac{Q_{n-2}(x)}{x^2Q_n(x)}y^n - \frac{1+x}{x^2(1-y)}.$$

\item Let $s(M)$ be the sum of all entries of the $m \times m$
matrix $M$ with nonnegative integral entries,  and let
$\bar{M}^{(k)}(x)$ be the $m \times m$ matrix $I-xM$ where the
$k$-th row of the matrix $I-xM$ is replaced by $(1,1, \cdots,1)$.
$f_k(x)$ be the determinant of the matrix  $\bar{M}^{(k)}(x)$.
Then the generating function $\eta(x)$ for the sequence
$\{s(M^n)\}_{n=0}^\infty$ is as follows (the proof is immediate
from the 4.7.2 Theorem of R. Stanley \cite{ST2}) :
$$ \eta(x)=\sum_{n=0}^\infty s(M^n)x^n = \frac{\sum_{i=1}^m
f_i(x)}{\det(I-xM)}.$$

\end{enumerate}
\end{rmk}
\bigskip
\subsection{Circular Graphs}\label{circGraph}

Before considering circular graphs, we first consider the linear
graphs $L_{k+1}$ with $k+1$ vertices, two of them
being end vertices $a$ and $b$, as follows below:\\
\xymatrix{a{\bullet} \ar@{-}[r] & {\bullet} \ar@{-}[r] &
{\bullet}\cdots{\bullet} \ar@{-}[r] & {\bullet}
 \ar@{-}[r] & {\bullet} \ar@{-}[r] & {\bullet}b\\
&&&&&}\\ If we fix $n_a =i$ and  $n_b=j$, then the number of all
 possible ways to give distributions corresponding to the
vertex set of the graph $L_{k+1}$ is same as the $(i,j)-$entry of
the matrix $B(n)^k$.

Next, we identify the leftmost vertex $a$ and the rightmost vertex
$b$, resulting in the circular graph $C_k$ with $k$ vertices.

Hence, the number of all possible ways to give distributions
corresponding to the vertex set of the graph $C_{k}$ (with
identified vertex having $n_a=n_b=i$ fixed) is same as the
$(i,i)$-entry of the matrix $B(n)^k$.

This implies that $WC_k(n)=tr(B(n)^k)$ for $k \ge 2$, and
$WC_1(n)=trace(I_{k+1})=n+1.$

\bigskip

\begin{lmma}[Theorem 4.7.2 and Corollary 4.7.3, Stanley, 1997
\cite{ST2}]\label{stanley1}

For any $m \times m$ matrix $M$, we have:
$$ \sum_{l=1}^\infty tr(M^l)x^l = - \frac{x
\frac{d}{dx} (\det(I-xM))}{\det(I-xM)}.$$
\end{lmma}
\bigskip

By the Lemma \ref{stanley1} above \quad
$$CF(n,x)=\sum_{k=1}^\infty
tr(B(n)^k)x^k = \lfloor \frac{n+1}{2} \rfloor x -
       \frac{x{Q_n}^{\prime}(x)}{Q_n(x)}$$
and
$$CG(x,y)=\sum_{k=1,n=0}^\infty WC_k(n)x^ky^n =
\frac{xy}{(1+y)(1-y)^2} - \sum_{n=0}^\infty
\frac{x{Q_n}^{\prime}(x)}{Q_n(x)}y^n.$$

We list several $\rho(C_i)$'s and $\rho(L_i)$'s for reference.

$\begin{array}{ll}
\rho(C_1)=\frac{1}{(1-x)^2}                           & \rho(L_1)=\frac{1}{(1-x)^2}        \\
\rho(C_2)=\frac{1}{(1-x)^3}                           & \rho(L_2)=\frac{1}{(1-x)^3}        \\
\rho(C_3)=\frac{1+x+x^2}{(1+x)(1-x)^4}                &  \rho(L_3)=\frac{1+x}{(1-x)^4}      \\
\rho(C_4)=\frac{1+2x+x^2}{(1-x)^5}                    &  \rho(L_4)=\frac{1+3x+x^2}{(1-x)^5} \\
\rho(C_5)=\frac{1+6x+11x^2+6x^3+x^4}{(1+x)(1-x)^6}
&\rho(L_5)=\frac{1+7x+7x^2+x^3}{(1-x)^6}
\end{array}$

\bigskip
For example, $WC_3(n)=\frac{1}{16} (4n^3+18n^2+28n+15+(-1)^n)\\
=\begin{cases}
        \frac{1}{8}(n+2)(2n^2+5n+4) & \text{if $n$ is even} \\
        \frac{1}{8}(n+1)(2n^2+7n+7) & \text{if $n$ is odd}
 \end{cases}$.

\subsection{Discrete Graphs}\label{disGraph}
In order to describe discrete graphs, complete graphs, star
graphs, complete bipartite graphs, etc, we need the notion of
 Eulerian numbers and Eulerian polynomials.

Let $p=p_1p_2 \cdots p_t$ be a $t$-permutation. We say that $i$ is
a descent of $p$ if $p_i > p_{i+1}.$ Let $A(t,k)$ be the number of
$t$-permutations with $k-1$ descents. The numbers $A(t,k)$ are
called the {\bf Eulerian numbers}, and $A_t(x)=\sum_{k=1}^{t}
A(t,k)x^k$ is called the {\bf Eulerian polynomial}.\\

\begin{fact}We list below
several known facts about Eulerian numbers and Eulerian
polynomials. [See B\'ona, 2004 and 2005 \cite{BO1, BO2} for
details about Eulerian numbers, and also Graham et. al., 1994
\cite{GKP}]\label{Facts}

\noindent
\begin{enumerate}
\item $A(t,k)=kA(t-1,k)+(t-k+1)A(t-1,k-1).$

\item $A(t,k)=A(t,t+1-k)$ and $\sum_{k=1}^{t} A(t,k)=t!.$

\item $x^t = \sum_{k=1}^{t} A(t,k)C(x-1+k,t),$ where
$C(t,k)=\frac{t!}{k!(t-k)!}.$

\item $g(x,y)=\sum_{t=0}^\infty \sum_{k=1}^{t} A(t,k)x^k
\frac{y^t}{t!}=\sum_{t=0}^\infty A_t(x) \frac{y^t}{t!}=
\frac{1-x}{1-xe^{y(1-x)}}.$

\item $ \frac{A_t(x)}{x(1-x)^{t+1}} = \frac{d}{dx} \{
\frac{A_{t-1}(x)}{(1-x)^t}\}.$

\item $A_t(x)=txA_{t-1}(x)+x(1-x)A_{t-1}^\prime(x).$

\item $A_t(x)=x(1-x)^t + \sum_{k=1}^t C(t,k)(1-x)^{t-k} A_k(x).$

\item $A(t,k)=\sum_{i=0}^k (-1)^i C(t+1,i)(k-i)^n.$
\end{enumerate}
\end{fact}
\bigskip

Let $D_t=(V,E)$ be a discrete graph, that is, a graph with $|V|=t$
and $E=\emptyset$. Then, clearly,
$$WD_t(n)=(n+1)^t, \quad \rho(D_t)=\frac{A_t(x)}{(1-x)^{t+1}},$$ and
$$F(x,y)=\sum_{t=0}^\infty \rho(D_t) \frac{y^t}{t!} =
\frac{1}{1-xe^y}.$$
\subsection{Complete Graphs}\label{compGraph}
Let $K_t$ be a complete graph of order $t$. Then,

\bigskip

\begin{thm}[Ju, 2006 \cite{JU}]

\medskip
$$ WK_t(n)=t \sum_{r=1}^{\lfloor{\frac{n+1}{2}}\rfloor} r^{t-1} +
(\lfloor{\frac{n+2}{2}}\rfloor)^t,$$ and
$$F(x,y):=\sum_{t=0}^\infty \sum_{n=0}^\infty WK_t(n)x^n \frac{y^t}{t!}
=(1+x+\frac{xy}{1-x}) \frac{e^y}{1-x^2e^y}. $$
\end{thm}
\bigskip

We define $$\alpha(y;x) = 1+x+\frac{xy}{1-x}, \qquad
 \beta(y;x)=\frac{e^y}{1-x^2 e^y}.$$

Using the facts given in Section \ref{disGraph} we can show that
the following holds:

\begin{equation*}
\begin{aligned}
\beta(y;x)&=\frac{e^y}{1-x^2 e^y}=(\sum_{t=0}^\infty
\frac{y^t}{t!})(\frac{1}{1-x^2} + \sum_{t=1}^\infty \frac{
A_t(x^2)}{(1-x^2)^{t+1}} \frac{y^t}{t!})\\
&=\frac{1}{1-x^2} + \sum_{t=1}^\infty
\frac{A_t(x^2)}{x^2(1-x^2)^{t+1}} \frac{y^t}{t!},
\end{aligned}
\end{equation*}
and
\begin{equation}\label{eq:2}
F(x,y)= \frac{1}{1-x}+\frac{y}{(1-x)^2}+\sum_{t=2}^\infty
      \frac{A_t(x^2)+txA_{t-1}(x^2)}{x^2(1-x)(1-x^2)^t}\frac{y^t}{t!}.
\end{equation}

For $t \geq 2$, the numerator of the summand in the equation
\eqref{eq:2} is
$$\frac{A_t(x^2)+txA_{t-1}(x^2)}{x^2}=(1+x)\frac{txA_{t-1}(x^2) +
x^2(1-x^2)A_{t-1}^\prime (x^2)}{x^2}.$$ Let
$$r_t(x)=\frac{txA_{t-1}(x^2) +
x^2(1-x^2)A_{t-1}^\prime (x^2)}{x^2},$$ for $t=2,3, \cdots$. By
Fact \ref{Facts} 2. (symmetric condition),
$$r_t(-1)=\sum_{k=1}^{t-1} (2k-t)A(t-1,k)=0,$$ for $t=2,3, \cdots$.
This implies that
$$\frac{A_t(x^2)+txA_{t-1}(x^2)}{x^2}=(1+x)^2 p_{2t-4}(x)$$ for some
polynomial $p_{2t-4}(x)$ in $x$ of degree $2t-4$.

\bigskip

\begin{thm}
For $U_t(x)=\sum_{n=0}^\infty WK_t(n)x^n,$
$$F(x,y)=\sum_{t=0}^\infty U_t(x) \frac{y^t}{t!}=\frac{1}{1-x}+
\frac{y}{(1-x)^2}+ \sum_{t=2}^\infty
\frac{p_{2t-4}(x)}{(1-x)^3(1-x^2)^{t-2}}
\frac{y^t}{t!},$$\mbox{where}
$$p_{2t-4}(x)=\frac{A_t(x^2)+txA_{t-1}(x^2)}{x^2(1+x)^2}$$ is a polynomial of
degree $2t-4$ for $t=2,3,4, \cdots. $
\end{thm}
\bigskip

Computations (using Maple) of $U_t(x)$ for $t=0,1,2,\cdots,8$ have
shown that:
\begin{equation*}
\begin{split}
U_0(x) & = \frac{1}{1-x},\\
U_1(x) & = \frac{1}{(1-x)^2},\\
U_2(x) & = \frac{1}{(1-x)^3},\\
U_3(x) & = \frac{1+x+x^2}{(1-x)^3 (1-x^2)},\\
U_4(x) & = \frac{1+2x+6x^2 + 2x^3 +x^4}{(1-x)^3 (1-x^2)^2}\\
U_5(x) & = \frac{1+3x+19x^2 + 14x^3 +19x^4 +3x^5 + x^6}{(1-x)^3 (1-x^2)^3}\\
U_6(x) & = \frac{1+4x+48x^2 + 56x^3 +142x^4 +56x^5 + 48x^6 + 4x^7 + x^8}{(1-x)^3 (1-x^2)^4}\\
U_7(x) & = \frac{1+5x+109x^2 + 176x^3 +730x^4 +478x^5 + 730x^6}{(1-x)^3 (1-x^2)^5}\\
& \qquad  + \frac{ 176x^7 + 109x^8 + 5x^9 +
x^{10}}{(1-x)^3(1-x^2)^5}\\
\end{split}
\end{equation*}

\begin{rmk}

\medskip

\noindent
\begin{enumerate}
\item The sequence $a_1,a_2, \cdots, a_n$ of positive real numbers
is called  {\it unimodal} if there exists an index $k$ such that
$1 \le k \le n$, and $a_1 \le a_2 \le \cdots \le a_k \ge a_{k+1}
\ge \cdots \ge a_n$. The same sequence is called {\it log-concave}
if $a_{k-1}a_{k+1} \le a_k^2$ holds for all indices $k.$ It is
well-known that a log-concave sequence is unimodal.(See
Bona\cite{BO2} for the proof.) The sequence $\{A(t,k\}_{k=1}^{t}$
of Eulerian numbers is log-concave, so it is unimodal for all $t$.
The coefficients of numerator in $U_t(x)$ is unimodal for
$t=1,2,3,4,6,$ but not always as shown in $U_5(x),U_7(x)$ above.


\item  The rational functions $U_t(x)$
 have  denominator $(1-x^2)$. This means that
$WK_t(n)$ is a (in fact, Ehrhart) {\bf quasi-polynomial} of a
certain polytope.(We will mention  this later in Section
\ref{genFun}.) Hence, its form depends on the parity of $n$.
 (See the next remark.)


\item The sequence $\{WK_t(n)\}_{t,n}$( for $t=0, 1, 2, \cdots, 5
)$ are provided below.
\begin{align*}
WK_0 (n) &= C(n,0) = 1\quad \left ( \mbox{empty graph} \right )\\
WK_1 (n) &= C(n+1,1)= n+1 \quad
\begin{pmatrix} \xymatrix{\bullet} \end{pmatrix}\\
WK_2 (n) &= C(n+2,2) = \frac{1}{2}( n^2 + 3n + 2 ) \quad
\begin{pmatrix} \xymatrix{{\bullet}\ar@{-}[r]&{\bullet}} \end{pmatrix}\\
WK_3 (n) &=  \frac{1}{16}(4n^3 + 18n^2 + 28n + 15 + (-1)^n ) \quad
\begin{pmatrix}
\xymatrix{     &{\bullet}\ar@{-}[dr] \\
              {\bullet} \ar@{-}[ur]&
                 & {\bullet} \ar@{-}[ll]
                 } \end{pmatrix}\\
&\text{( Remarks 4 of Section 6 of Bona and Ju\cite{BJ}) }\\
WK_4 (n) &=  \frac{1}{16}(2n^4 + 12n^3 + 28n^2 + 30 n + 13\\
   & + (-1)^n (2n + 3) )
   \quad \begin{pmatrix}\xymatrix{ {\bullet}
\ar@{-}[d]\ar@{-}[dr] \ar@{-}[r]
   & {\bullet} \ar@{-}[d]\ar@{-}[dl]\\
   {\bullet}  & {\bullet}  \ar@{-}[l]
   } \end{pmatrix}\\
WK_5 (n) &=  \frac{1}{192}(12n^5 + 90n^4 +280n^3 +450 n^2 +
374n + 129 \\
   & + (-1)^n (30n^2 + 90n +63))
\end{align*}
\end{enumerate}
\end{rmk}
\subsection{Star Graphs}\label{starGraph}
The {\em Star Graph}
 $S_t$ of order $t$ is a tree with $t+1$ vertices, one
of them of degree $t$ and all others of degree $1$.\\
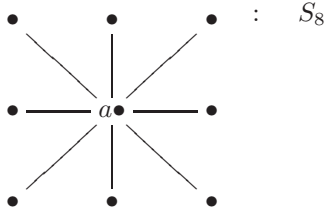
\xymatrix{{\bullet}
 \ar@{-}[dr] &{\bullet} \ar@{-}[d]  &{\bullet} \ar@{-}[dl]\\
          {\bullet} \ar@{-}[r]  &{a\bullet} \ar@{-}[r]  &{\bullet}\\
          {\bullet} \ar@{-}[ur] &{\bullet} \ar@{-}[u]  &{\bullet}
          \ar@{-}[ul]}
\quad : \quad $S_8$

If we let the hub vertex $a$ have  value $i$, then the rest of all
the vertices must have values in $[n-i]$. So
\begin{equation*}
\begin{split}
WS_t(n)&=\sum_{i=0}^n (n+1-i)^t =\sum_{k=1}^{n+1} k^t
=\sum_{k=1}^{n+1} [ \sum_{i=1}^{t}A(t,i)C(k-1+i,t)]\\
&=\sum_{i=1}^{t} A(t,i)C(n+1+i,t+1)
\end{split}
\end{equation*}
by Fact \ref{Facts} 3. in Section \ref{disGraph}.

 Now $$\rho(S_t)=\sum_{n=0}^\infty WS_t(n)x^n =
\sum_{i=1}^{t} \frac{x^{t-i}}{(1-x)^{t+2}}
A(t,i)=\frac{A_t(x)}{x(1-x)^{t+2}},$$ and

$$F(x,y)=\sum_{t=0}^\infty \rho(S_t) \frac{y^t}{t!} =
\sum_{t=0}^\infty \frac{A_t(x)}{x(1-x)^{t+2}} \ \frac{y^t}{t!}
=\frac{1}{x(1-x)(1-xe^y)},$$ by Fact \ref{Facts} 4.

\subsection{Cubic Graphs}\label{cubeGraph}
 The {\it Cubic Graph} $C$ is defined as
follows: 
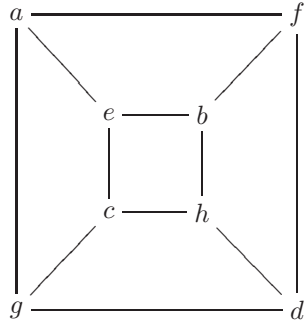
\xymatrix{ a \ar@{-}[rrr] \ar@{-}[ddd] \ar@{-}[dr] &&&f
\ar@{-}[ddd] \\
& e \ar@{-}[r] \ar@{-}[d] & b \ar@{-}[ur] \ar@{-}[d] & \\
& c \ar@{-}[r] \ar@{-}[dl] & h \ar@{-}[dr]  & \\
g \ar@{-}[rrr]  &&&d }\\

Given four numbers $a,b,c,d \in [n], e \in [n-max(a,b,c)],\\ f \in
[n-max(a,b,d)], g \in [n-max(a,c,d)],$ and  $h \in
[n-max(b,c,d)].$

\medskip
Hence
\begin{align*}
WC(n)&=\sum_{a,b,c,d=0}^n (n+1-max(a,b,c))(n+1-max(a,b,d))\\
     &(n+1-max(a,c,d))(n+1-max(b,c,d))\\
     &=\sum_{d=0}^n \sum_{max(a,b,c)=0}^n(n+1-max(a,b,c))(n+1-max(a,b,d))\\
     &(n+1-max(a,c,d))(n+1-max(b,c,d))\\
     &=\sum_{d=0}^n (\sum_{k=0}^d + \sum_{k=d+1}^n
)\sum_{max(a,b,c)=k}(n+1-max(a,b,c))\\
     &(n+1-max(a,b,d))(n+1-max(a,c,d))(n+1-max(b,c,d))\\
     &=\sum_{d=0}^n \sum_{k=0}^d \sum_{max(a,b,c)=k}(n+1-d)^3(n+1-k)\\
     &+\sum_{d=0}^{n-1} \sum_{k=d+1}^n \sum_{max(a,b,c)=k} (n+1-k)\\
     &(n+1-max(a,b))(n+1-max(a,c))(n+1-max(b,c))\\
     &:=[I]+[II].
\end{align*}

\medskip
$$[I]=
\frac{1}{3360}(n+1)(n+2)(n+3)(15n^5+138n^4+533n^3+1074n^2+1180n+560).$$
By the Inclusion Exclusion Principle,
$$\sum_{max(a,b,c)=k}(n+1-k)(n+1-max(a,b))(n+1-max(a,c))(n+1-max(b,c))$$
$$=3\sum_{l=0}^k (n+1-k)^3(n+1-l)((l+1)^2-l^2) -3\sum_{l=0}^k
(n+1-k)^4 + (n+1-k)^4.$$

$$[II]=\frac{1}{3360}
n(n+1)(n+2)(45n^5+357n^4+1177n^3+1971n^2+1638n+412).$$

\begin{align*}
WC(n)&=[I]+[II]\\
     &=C(n+8,8)+26C(n+7,8)+175C(n+6,8)+316C(n+5,8)\\
     &+175C(n+4,8)+26C(n+3,8)+C(n+2,8).\end{align*}

$$\rho(C)=\frac{1+26x+175x^2+316x^3+175x^4+26x^5+x^6}{(1-x)^9}.$$

Let $HC_d$ be a $d$-dimensional Hypercubic Graph. Then

\begin{align*}
\rho(HC_0)&= \frac{1}{(1-x)^2}\\
\rho(HC_1)&= \frac{1}{(1-x)^3}\\
\rho(HC_2)&= \frac{1+2x+x^2}{(1-x)^5}\\
\rho(HC_3)&=\frac{1+26x+175x^2+316x^3+175x^4+26x^5+x^6}{(1-x)^9}\\
\rho(HC_d)&=\frac{P_{2^d -2}(x)}{(1-x)^{1+2^d}} (d \ge 1)
\end{align*} (We conjecture the last one(general case)!)
where $P_{2^d -2}(x)$ is a symmetric polynomial of degree $2^d
-2$.

\subsection{Octahedral Graphs}\label{OctaGraph}

In this section we focus on octahedral graphs. An octahedral graph
is the Platonic graph with six nodes and 12 edges having the
connectivity of the octahedron. Let us consider the following
graph.

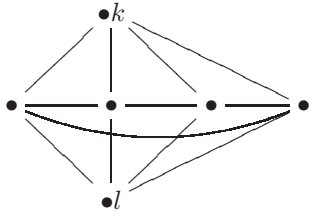
\xymatrix{ & {\bullet}k \ar@{-}[dl] \ar@{-}[d] \ar@{-}[dr]
\ar@{-}[drr] & & \\
{\bullet} \ar@{-}[r] \ar@{-}[dr] \ar@{-}@/_1pc/[rrr] & {\bullet}
\ar@{-}[r] \ar@{-}[d] & {\bullet} \ar@{-}[r] \ar@{-}[dl] &
{\bullet} \ar@{-}[dll]
\\ &{\bullet}l && }

Let $OH$ be an octahedral graph given  in the figure above. Let a
top vertex have value $k$ and a bottom vertex $l$ as in the
figure. Then rest of them have values in $[n-m]$, where
$m=max(k,l)$. Let $B(n,m)=(b(n,m)_{ij})$ be an $(n+1)\times(n+1)$
matrix for which $b(n,m)_{ij}=0$ if $max(i+j,max(i,j)+m)>n$ and
$b(n,m)_{ij}=1$ otherwise. For example,

\medskip

$B(4,0) =\begin{bmatrix}
1 & 1 & 1 & 1 & 1 \\
1 & 1 & 1 & 1 & 0 \\
1 & 1 & 1 & 0 & 0 \\
1 & 1 & 0 & 0 & 0 \\
1 & 0 & 0 & 0 & 0
\end{bmatrix},
B(4,1) =\begin{bmatrix}
1 & 1 & 1 & 1 & 0 \\
1 & 1 & 1 & 1 & 0 \\
1 & 1 & 1 & 0 & 0 \\
1 & 1 & 0 & 0 & 0 \\
0 & 0 & 0 & 0 & 0
\end{bmatrix},
B(4,2) =\begin{bmatrix}
1 & 1 & 1 & 0 & 0 \\
1 & 1 & 1 & 0 & 0 \\
1 & 1 & 1 & 0 & 0 \\
0 & 0 & 0 & 0 & 0 \\
0 & 0 & 0 & 0 & 0
\end{bmatrix}$,\\

\medskip

$B(4,3) =\begin{bmatrix}
1 & 1 & 0 & 0 & 0 \\
1 & 1 & 0 & 0 & 0 \\
0 & 0 & 0 & 0 & 0 \\
0 & 0 & 0 & 0 & 0 \\
0 & 0 & 0 & 0 & 0
\end{bmatrix},$
$B(4,4) =\begin{bmatrix}
1 & 0 & 0 & 0 & 0 \\
0 & 0 & 0 & 0 & 0 \\
0 & 0 & 0 & 0 & 0 \\
0 & 0 & 0 & 0 & 0 \\
0 & 0 & 0 & 0 & 0
\end{bmatrix}.$

\begin{align*}
WOH(n)&=\sum_{m=0}^n ((m+1)^2-m^2) trace(B(n,m)^4)\\
WOH(2k)&=\sum_{m=0}^{k-1} (2m+1) trace(B(2k,m)^4)
+\sum_{m=k}^{2k} (2m+1)(2k+1-m)^4\\
&=:(I)+(II)\\
WOH(2k+1)&=\sum_{m=0}^{k} (2m+1) trace(B(2k+1,m)^4)
 + \sum_{m=k+1}^{2k+1} (2m+1)(2k+2-m)^4\\
 &=:(III)+(IV).
 \end{align*}

In order to compute the quantities (I) and (III) above, we need to
define the following.  
Let $D(n)$ be an $(n+1)\times (n+1)$ matrix all of whose entries
are 1, and let $S(n,m):=(s_{ij})$ be an $(n+1)\times (n+1)$ matrix
such that $s_{ij}=1$ if $i+j>2n-m$ and $0$ otherwise.

\medskip

For example, $$D(4)-S(4,2)=
 \begin{bmatrix}
1 & 1 & 1 & 1 & 1 \\
1 & 1 & 1 & 1 & 1 \\
1 & 1 & 1 & 1 & 1 \\
1 & 1 & 1 & 1 & 0 \\
1 & 1 & 1 & 0 & 0
\end{bmatrix}
=\begin{bmatrix}
1 & 1 & 1 & 1 & 1 \\
1 & 1 & 1 & 1 & 1 \\
1 & 1 & 1 & 1 & 1 \\
1 & 1 & 1 & 1 & 1 \\
1 & 1 & 1 & 1 & 1
\end{bmatrix}
-\begin{bmatrix}
0 & 0 & 0 & 0 & 0 \\
0 & 0 & 0 & 0 & 0 \\
0 & 0 & 0 & 0 & 0 \\
0 & 0 & 0 & 0 & 1 \\
0 & 0 & 0 & 1 & 1
\end{bmatrix}.$$

\begin{align*}
(D(r)-S(r,m))^4 &= (D(r)^2 -D(r)S(r,m) - S(r,m)D(r,m) +
S(r,m)^2)^2 \\
&=:(M-DS-SD+K)^2 \\ &= M^2 - MDS - MSD + MK - DSM +DSDS + DSSD - DSK \\
&- SDM + SDDS + SDSD - SDK + KM - KDS - KSD + K^2.
\end{align*}

Let us define
\begin{align*}
p(r,m)&:=trace((D(r-1)-S(r-1,m))^4)\\
&=r^4-4C(m+1,2)r^2+4[C(m+1,3)+C(m+2,3)]r\\
&-4C(m+2,4)-C(m+1,2).
\end{align*}

Then we have
\begin{align*}
WOH(2k)&=\sum_{m=0}^{k-1} (2m+1)p(2k+1-m,2(k-m))\\
 &+\sum_{m=k}^{2k} (2m+1)(2k+1-m)^4\\
&=\frac{1}{10}(k+1)(2k^2+2k+1)(12k^3+30k^2+27k+10),\\
WOH(2k+1)&=\sum_{m=0}^{k} (2m+1)p(2k+2-m,2k+1-2m)\\
&+\sum_{m=k+1}^{2k+1} (2m+1)(2k+2-m)^4\\
&=\frac{1}{10}(k+1)(2k^2+6k+5)(12k^3+42k^2+51k+20).
\end{align*}

Finally, we get
\begin{align*}
WOH(n)&=\sum_{m=0}^{\lfloor \frac{n-1}{2} \rfloor}
(2m+1)p(n+1-m,n-2m)\\&+ \sum_{m=\lfloor \frac{n+1}{2} \rfloor}^n (2m+1)(n+1-m)^4\\
&=\frac{1}{160}(6n^6+54n^5+210n^4+450n^3+559n^2+381n+115\\&+(-1)^n
(10n^3+45n^2+75n+45))
\end{align*}
and
$$\rho(OH)=
\frac{1+7x+48x^2+89x^3+142x^4+89x^5+48x^6+7x^7+x^8}{(1+x)^4
(1-x)^7}$$

\subsection{Complete Bipartite Graphs}\label{comBiGraph}

Let $K_{p,q}=(V,E)$ be a complete bipartite graph of order
$(p,q)$. So, $V=X \cup Y, E=\{ab=ba| a \in X, b\in Y
\},|X|=p,|Y|=q.$ We also let $WK_{p,q,\alpha}$ be a weighted graph
with a distribution $\alpha$ of the form
$\alpha=((n_1,n_2,\cdots,n_p),(m_1,m_2,\cdots,m_q))$ and
$r=max(n_1,n_2,\cdots,n_p).$ Then $m_j \in [n-r]
(j=1,2,\cdots,q)$.

\begin{align*}
WK_{p,q}(n)&=\sum_{n_1,n_2,\cdots,n_p =0}^n
(n+1-max(n_1,n_2,\cdots,n_p))^q\\
&=\sum_{k=0}^n ((k+1)^p-k^p)(n+1-k)^q
\end{align*}
If we let $\bar{A}(t,k)=A(t,k+1)$ and
$\bar{A}_t(x)=\sum_{k=0}^{t-1} \bar{A}(t,k)x^k.$ Then
$$\rho(K_{p,q})=\sum_{n=0}^\infty
|WK_{p,q}(n)|x^n=\frac{\bar{A}_p(x) \bar{A}_q(x)}{(1-x)^{p+q+1}}$$

For example,
\begin{align*}
\bar{A}_3(x)&=\bar{A}(3,0)+\bar{A}(3,1)x+\bar{A}(3,2)x^2\\&=1+4x+x^2, \\
\bar{A}_4(x)&=\bar{A}(4,0)+\bar{A}(4,1)x+\bar{A}(4,2)x^2+\bar{A}(4,3)x^3
\\&=1+11x+11x^2+x^3,
\end{align*}
and
$$\rho(K_{3,4})=\rho(K_{4,3})=
\frac{(1+4x+x^2)(1+11x+11x^2+x^3)}{(1-x)^8}.$$

\begin{rmk}

\medskip

\noindent
\begin{enumerate}
\item We have $H(x,y,z):= \sum_{p,q=0}^\infty \rho(K_{p,q}) \frac{y^p}{p!}
\frac{z^q}{q!} = \frac{1-x}{(1-xe^y)(1-xe^z)}$

\item We have $\rho(K_{p,p})=\frac{\bar{A}_p(x)^2}{(1-x)^{2p+1}}$.
\item The graph $K_{p,q}$ has no cycles of odd length. Hence
generating functions are not quasi-polynomials (as we will show
later).

\item The degree of the numerator in $\rho(K_{p,q})$ is $p+q-2$
and the numerator is symmetric.
\end{enumerate}
\end{rmk}

\section{Rational Polytopes and Rational Generating Functions}\label{genFun}

Readers unfamiliar with this topic may wish to consult the book of
R. Stanley(\cite{ST2}, Chapter 4). A {\bf quasi (or
pseudo)-polynomial} of degree $d$ with a quasi-period $N$ is a
function $f:\N \rightarrow \C$ of the form $f(n)=\sum_{i=0}^d
c_i(n)n^i$ where the coefficients are periodic functions of a
common period $N$ and the leading coefficient $c_d(n)$ is not
identically zero.
\medskip

\begin{ex}\label{ex3:1}
The number $f(n)$ of unit squares in the region bounded by $x=0,
x=n(\ge 1), y=0, y=\frac{4}{3}x$ :

\medskip

$f(n)=\begin{cases}
   \frac{4}{3}n-2 &\text{if}\ n \equiv 0(mod3)\\
   \frac{4}{3}n-\frac{4}{3} &\text{if}\ n \equiv 1(mod3)\\
   \frac{4}{3}n-\frac{5}{3} &\text{if}\ n \equiv 2(mod3)
   \end{cases}$.\\ $c_1(n)=\frac{4}{3}$ is constant(1-periodic),
 but $c_0(n)$ is 3-periodic.
\medskip

The generating function associated with the sequence $\{
f(n)\}_{n=0}^{\infty}$ is as follows: $$g(x)=\sum_{n=0}^\infty
f(n)x^n = \frac{x^2(1+x+2x^2)}{(1-x)(1-x^3)}. $$
\end{ex}

Here is another nontrivial example from the page 220 (Exercises
13.9 and 13.10) of Pach and Agarwal \cite{PA}.

\begin{ex}\label{ex3:2}
Given a set $P$ of $n$ points in the plane, for any $p \in P,$ let
$\mu_P(p)$ denote the number of farthest neighbors of $p$, that
is,
$$\mu_P(p)=|\{q \in P||p-q|=max_{r \in P} |p-r| \}|.$$

Let $\mu(n)=max_{|P|=n} \sum_{p \in P} \mu_P(p)$. It turned out
(see Csizmadia \cite{CG} or Avis et al \cite{AEP}) that if $n$ is
sufficiently large, then
$$\mu(n)=\frac{n^2}{4}+\frac{3n}{2}+ \begin{cases}
       3           &\text{if } n \equiv 0 (mod2),\\
       \frac{9}{4} &\text{if } n \equiv 1 (mod4),\\
       \frac{13}{4} &\text{if } n \equiv 3 (mod4).\\
       \end{cases}$$
The generating function associated with the sequence $\{
\mu(n)\}_{n=0}^{\infty}$ is as follows:
$$g(x)=\sum_{n=0}^\infty \mu(n)x^n =
\frac{3-3x+4x^2-x^3-3x^4+2x^5}{(1-x)^2(1-x^4)}.
$$
\end{ex}
\smallskip

If denominator in the reduced generating function associated with
a certain sequence $\{ f(n)\}_{n=0}^{\infty}$ has a factor $1-x^N$
then $f(n)$ has a quasi-period $N$ as in Example \ref{ex3:1}
($1-x^3$) and Example \ref{ex3:2} ($1-x^4$) \cite{ST2}.

\medskip

Suppose we have a finite set of points $\{t_1,t_2, \cdots,t_d \}$
in $\R^m$. The {\bf convex hull} of the set $\{t_1,t_2, \cdots,t_d
\}$ is the set of all convex combinations of the given points,
i.e. $\{x \in \R^m: x = \sum_{i=1}^d \lambda_i t_i, \, \lambda
\geq 0, \, \sum_{i=1}^d \lambda_i = 1\}$. By a {\bf convex
polytope}, or simply a {\bf polytope}, we mean a set which is the
convex hull of a non-empty finite set $\{t_1,t_2, \cdots,t_d \}.$
An {\bf affine combination} of points $t_1, t_2, \cdots, t_k$ from
$\R^m$ is a linear combination $\lambda_1 t_1+\lambda_2 t_2+
\cdots + \lambda_k t_k$, where $\lambda_1 + \cdots + \lambda_k =1$
and $\lambda_i \in \R$ for $i = 1, \ldots , k$. A $k$-family
$(t_1,t_2,\cdots,t_k)$ of points from $\R^m$ is said to be {\bf
affinely independent} if a linear combination $\lambda_1
t_1+\lambda_2 t_2+ \cdots + \lambda_k t_k$ with $\lambda_1 +
\cdots + \lambda_k =0$ can only have the value $0$ when
$\lambda_1=\lambda_2=\cdots=\lambda_k=0.$ A polytope ${\mathcal
P}$ with the property that there exists an affinely independent
family $(t_1,t_2, \cdots,t_d )$ such that ${\mathcal P}$ is a
convex hull of $\{t_1,t_2, \cdots,t_d \}$ is called a {\bf
simplex}. (Refer Barvinok \cite{BA}, Br{\o}ndsted \cite{BR},
Miller and Sturmfels \cite{MS}, Stanley \cite{ST1} or Schrijver
\cite{schrijver} for details on polytopes.)

For a graph $G=(V,E) \in {\mathcal G}$, let $m=|V|$ and
$${\mathcal P}(G) = \{ (x_1,x_2,\cdots,x_m) \in
\R_{+}^m | x_i+x_j \le 1 \text{ for all edge } v_iv_j \in E \}.$$
Note that for every simple graph $G$ ${\mathcal P}(G)$ is a
polytope which is contained in the $m$-dimensional unit hypercube.
If all  of  the coordinates of the vertices of the polytope are
integers, then we call such a polytope an {\bf integer polytope}
(or an {\bf integral polytope}, a {\bf lattice polytope}). If all
of the coordinates of vertices of the polytope are rational
numbers, then the associated polytope is called a {\bf rational
polytope}.

Note that every ${\mathcal P}$ is homeomorphic to a ball ${\bf
B}^d$, for some $d$. This $d$ is the dimension $dim({\mathcal P})$
of the polytope ${\mathcal P}$. We denote the {\bf boundary}
(resp. {\bf interior}) of ${\mathcal P}$ by ${\partial {\mathcal
P}}$(resp. ${\bar {\mathcal P}}$). ${\bf \alpha} \in {\mathcal P}$
is a {\bf vertex} of ${\mathcal P}$ if there exists a closed
affine half-space ${\mathcal H}$ such that
${\mathcal P} \cap {\mathcal H} = \{\alpha \}$\\


If ${\mathcal P} \in \R^m$ is a rational convex polytope and if $n
\in \Z_+$, then we define $i({\mathcal P},n):=|n{\mathcal P} \cap
\Z^m|$ and
${\bar i}({\mathcal P},n):=|n{\bar {\mathcal P}} \cap \Z^m|$.\\
This is called the {\bf Ehrhart quasi-polynomial} of ${\mathcal
P}$ and ${\bar {\mathcal P}}$ respectively. (We will show later
that this is a quasi-polynomial.) If ${\mathcal P}$ is an integer
polytope then this is simply called the {\bf Ehrhart polynomial}
of ${\mathcal P}$ and ${\bar {\mathcal P}}$ respectively.

If $\beta \subset \Q^m$, then define  $den{\bf \beta}$ (the
denominator of ${\bf \beta}$) as the least integer $q \in \Z_+$
such that $q\beta \subset \Z^m$.

\bigskip

\begin{thm}[Stanley, 1997 \cite{ST2}, pp237-238]\label{Stanley2}
If ${\mathcal P}$ is a rational convex polytope of dimension $d$
in $\R^m$, and
$$F({\mathcal P},x):=1+\sum_{n=1}^\infty i({\mathcal P},n) x^n.$$
Then $F({\mathcal P},x)$ is a rational function
$\frac{P(x)}{Q(x)}$ of $x$, where $deg(P(x))\le d$ and $Q(x)$ {\bf
can be} written as $\prod_{{\bf \beta} \in V} (1-x^{den {\bf
\beta}})$. If $F({\mathcal P},x)$ is written in lowest terms, then
$x=1$ is a pole of order $d+1$, and no value of $x$ is a pole of
order $> d+1$. We also have the following (reciprocity for Ehrhart
quasi-polynomial):
$$\bar{i}({\mathcal P},n)=(-1)^d i({\mathcal P},-n).$$ 
\end{thm}
\bigskip

Theorem \ref{Stanley2} says that $i({\mathcal P},n)$ is a {\bf
quasi-polynomial} with the correct value $i({\mathcal P},0)=1$,
and $D(x)=\prod_{{\bf \beta} \in V} (1-x^{den {\bf \beta}})$ is
not in general the least denominator of $F({\mathcal P},x)$.
However, the least denominator has a factor $(1-x)^{d+1}$ but not
$(1-x)^{d+2}$, while $D(x)$ has a factor $(1-x)^{|V|}$. See the
following example (Stanley \cite{ST2}).

\medskip

\begin{ex}
Let ${\mathcal P}$ be the convex hull of the vertices set\\
$V^*= \{(0,0,0),(1,0,0),(0,1,0),(1,1,0),(1/2,0,1/2)\}$. Then
$$D(x)=\prod_{{\bf \beta} \in V^*} (1-x^{den {\bf
\beta}})=(1-x)^5(1+x),$$ but $F({\mathcal P},x)=\frac{1}{(1-x)^4}$
\end{ex}

\medskip

In order to prove our main theorem we need some information about
the denominator of the vertex set for the polytope ${\mathcal
P}(G)$.

\medskip

\begin{lmma}\label{Lemma:1}
Let $G=(V,E)$ with $|V|=m$ be a simple bipartite graph and $V^*$
the set of vertices of the polytope ${\mathcal P}(G).$ Then
$den(V^*)=\{1\}$. Otherwise (equivalently, if the simple graph $G$
is not bi-colorable, or if it has a cycle of odd length), then
$den(V^*)=\{1,2\}$.
\end{lmma}


\begin{proof}
If $G = (V, E)$ is bipartite.  Then the defining matrix for the
polytope ${\mathcal P}(G)$ is totally unimodular (Example 1 on p.
273, \cite{schrijver}). Thus, $den(V^*)=\{1\}$.

If $G = (V, E)$ is not bipartite. then, $G$ contains at least a
cycle with odd length. Consider a cycle with odd length as a
subgraph and a submatrix in the defining matrix for the polytope
${\mathcal P}(G)$.

Let $A$ be the defining matrix for $G$, and suppose the polytope
${\mathcal P}(G)$ is defined by the system $\{x \in \R^m: Ax \leq
b, \, x \geq 0\}$ where $A$ is a $s \times m$ matrix.  Then we can
rewrite this system as $\{x \in
\R^m: A'x \leq b'\}$ where $A' = \left(\begin{array}{c}A\\-I_m\\
\end{array} \right)$ where $I_m$ is an identity matrix of order
$m$ and $b' = \left(\begin{array}{c}b\\0
\\ \end{array} \right)$ where $0$ is a $m$ dimensional vector with all zeros.

Then note that any vertex of the polytope ${\mathcal P}(G)$ is
determined by the unique solution of a subsystem
\[
A'' x = b''
\]
where $A''$ is a $m \times m$ minor of $A'$ with $\det(A'') \neq
0$ and $b''$ is a $m$ sub-vector of $b'$ (Theorem 8.4
\cite{schrijver}). Thus all we have to show is that $\det(A'') =
\{\pm 1, \, \pm 2\}$ if  $\det(A'') \neq 0$. Then, there are three
cases we have to consider.
\begin{enumerate}
\item\label{s1} All rows of $A''$ are from $I_m$. \item\label{s2}
All rows of $A''$ are from $A$.
\begin{enumerate}
\item \label{case1} $A''$ contains a cycle (or cycles) with odd
length. \item \label{case2} otherwise.
\end{enumerate}
\item\label{s3} $k$ rows of $A''$ are from $A$ and $m - k$ rows
are from $I_m$.
\begin{enumerate}
\item \label{case21} $A''$  contains a cycle (or cycles) with
smaller odd length $k$. \item \label{case22} otherwise.
\end{enumerate}
\end{enumerate}
If $A'' = I_m$, then clearly this defines the origin. For Case
\ref{case2} and Case \ref{case22}, $A''$ does not contain a cycle
with odd length.  Thus we are done. For Case \ref{case1} and Case
\ref{case21}, we have to show that the matrix for each subsystem
of the system $A'' x = b''$ defining a cycle with its length $k$,
where $k$ is an odd positive integer with $k \geq 3$, has its
determinant $\pm 2$. Let $\bar A$ be the defining matrix for a
cycle with odd length. Note that $\bar A$ is a $k \times k$ matrix
and note that after permuting columns and rows, $\bar A$ forms
such that:

\[
\bar A =  \left( \begin{array}{cccccc}
1 & 1 & 0 & \cdots & 0 & 0\\
0 & 1 & 1 & \cdots & 0 & 0\\
\vdots & \vdots & \vdots & \vdots & \vdots & \vdots \\
0 & 0 & 0 & \cdots & 1 & 1\\
1 & 0 & 0 & \cdots & 0 & 1\\
\end{array}
\right)
\]

Note that $\bar A$ is an upper triangular matrix except for the
$1$ in the first position of the $k$th row.  So the determinant of
$\bar A$ is clearly $2$. Thus, $den(V^*)=\{1,2\}$.
\end{proof}


\medskip

The next theorem is the conclusion of all that we discussed. Its
proof follows immediately from Theorem \ref{Stanley2} and Lemma
\ref{Lemma:1}.

\medskip
\begin{thm}\label{thm:3:6}
If a simple graph $G=(V,E)$ with $|V|=m$ is bipartite
(equivalently, $\chi(G) \le 2$), then $WG(n)=i({\mathcal P}(G),n)$
is an Ehrhart polynomial and $$\rho(G)=\frac{P(x)}{(1-x)^{m+1}},$$
where $P(x)$ is a symmetric polynomial of degree $\le m$. That is,
if $WG(n)$ is an Ehrhart quasi-polynomial in $n$ of a quasi-period
$\neq 1$, then the quasi-period is 2, $G$ contains a cycle of odd
length (hence $\chi(G) > 2$) and
$$\rho(G)=\frac{P(x)}{(1-x^2)^k (1-x)^{m+1-k}},$$
where $P(x)$ is a polynomial of degree $\le m+k$ and $k$ is a
nonnegative integer.
\end{thm}

\medskip

\begin{cor}
If a simple graph $G=(V,E)$ with $|V|=m$ is either a tree, a
circular graph of even length, a discrete graph, a hypercubic
graph, a complete bipartite graph, or a grid graph, then $\rho(G)=
\frac{p(x)}{(1-x)^{m+1}}$, where $p(x)$ is a symmetric monic
polynomial of degree at most $m$.
\end{cor}
\begin{proof}
All of those graphs are bipartite graphs. Hence, by Theorem
\ref{thm:3:6}, the result follows.
\end{proof}
\medskip

\begin{ex}\label{eg3:8}
Consider the Circular Graph $C_3$ (a cycle of odd length 3).

$\begin{aligned}
&{\mathcal P}(C_3)=\{(r,s,t) \in \R^3 | r,s,t\ge 0, r+s,s+t,t+r \le 1\}\\
&V^* = \{(0,0,0),(1,0,0),(0,1,0),(0,0,1),(\frac{1}{2},\frac{1}{2},\frac{1}{2})\}\\
&m=dim({\mathcal P}(C_3))=3\\
&WC_3(n)=i({\mathcal P}(C_3),n)=\frac{1}{16}(4n^3+18n^2+28n+15+(-1)^n)\\
&\rho(C_3)=F({\mathcal P}(C_3),x)=1+\sum_{n=1}^\infty i({\mathcal
P}(C_3),n) x^n = \frac{1+x+x^2}{(1-x^2)(1-x)^3}
=\frac{P(x)}{(1-x^2)^1(1-x)^{3+1-1}}\\
&D(x)=\prod_{\beta \in V^*} (1-x^{den\beta})=(1-x^2)(1-x)^4.
\end{aligned}$
\end{ex}
\medskip

\begin{ex}\label{eg3:9}
Consider the Circular Graph $C_4$ (a cycle of even length 4).

$\begin{aligned}
&{\mathcal P}(C_4)=\{(r,s,t,u) \in \R^4 | r,s,t,u\ge 0, r+s,s+t,t+u,u+r \le 1\},\\
&V^* =\{(0,0,0,0),(1,0,0,0),(0,1,0,0),(0,0,1,0),(0,0,0,1),(0,1,0,1),(1,0,1,0)\},\\
&m=dim({\mathcal P}(C_4))=4,\\
&WC_4(n)=C(n+2,4)+2C(n+3,4)+c(n+4,4)=1+\frac{5}{2}n+\frac{7}{3}n^2+n^3+\frac{1}{6}n^4,\\
&\rho(C_4)=1+\sum_{n=1}^\infty i({\mathcal P}(C_4),n) x^n =
\frac{1+2x+x^2}{(1-x)^5}
=\frac{P(x)}{(1-x)^{4+1}},\\
&D(x)=\prod_{\beta \in V^*} (1-x^{den \beta})=(1-x)^7.
\end{aligned}$
\end{ex}
\medskip


\begin{ex}\label{eg3:10}
Consider the simple graph $G=(V,E)$, where $V=\{1,2,3,4,5\}$ and
$E=\{12,23,34,45,15,13\}.$ Hence, the associated polytope and its
vertices are (we used {\sffamily LattE} \cite{DHHHTY, latte1} and
{\sffamily CDD} \cite{FU} for computational experimentation):

$\begin{aligned}
&{\mathcal P}(G)= \{(r,s,t,u,v) \in \R^5 | r,s,t,u,v \ge 0, r+s,s+t,t+u,u+v,v+r,r+t \le 1 \}.\\
&V^*=\{(0,0,0,0,0),(1,0,0,0,0),(0,1,0,0,0),(0,0,1,0,0),(0,0,0,1,0),\\
&(0,0,0,0,1),(1,0,0,1,0),(0,1,0,1,0),(0,1,0,0,1),(0,0,1,0,1),\\
&(1/2,1/2,1/2,0,0),(1/2,1/2,1/2,1/2,0),(1/2,1/2,1/2,0,1/2),\\
&(1/2,1/2,1/2,1/2,1/2)\},\\
&dim({\mathcal P}(G))=5,\\
&WG(n)=\frac{121}{128}+\frac{535}{192}n+\frac{219}{64}n^2+\frac{13}{6}n^3+\frac{45}{64}n^4+\frac{3}{32}n^5
+(-1)^n(\frac{7}{128}+\frac{3}{64}n+\frac{1}{64}n^2),\\
&\rho(G)= \frac{1+7x+22x^2+30x^3+22x^4+7x^5+x^6}{(1-x)^3(1-x^2)^3},\\
&D(x)=\prod_{\beta \in V^*} (1-x^{den \beta})=(1-x)^{10}(1-x^2)^4.
\end{aligned}$
\end{ex}
\bigskip

$C_3$ in Example \ref{eg3:8} has a cycle of odd length (3) and
$C_4$ in Example \ref{eg3:9} has a cycle of even length (4).
However, the graph $G$ in Example \ref{eg3:10} has three cycles of
length 3, 4 and 5, two of them odd.

\medskip

\begin{rmk}
If a simple graph $G=(V,E)$ with $|V|=m$ is bipartite and
$\rho(G)=\frac{P(x)}{(1-x)^{m+1}},$ then coefficient vector
$h(G)=(h_0 ,h_1, \cdots , h_m)$ of a symmetric polynomial
$P(x)=h_0 + h_1x + \cdots + h_mx^m$ satisfies
the following properties(see Stanley \cite{ST3}):\\
(1) $h_0=1.$\\
(2) $h_m=(-1)^m WG(-1).$\\
(3) $min\{j \ge 0 | WG(-1)=WG(-2)=\cdots=WG(-(m-j))=0\}.$\\
    $=max\{i|h_i \ne 0 \}$\\
    For example, $WC_4(-1)=WC_4(-2)=0$ and $h(C_4)=(1,2,1,0,0).$\\
(4) $WG(-n-k)=(-1)^m WG(n)$ for all $n$ if and only if
$h_i=h_{m+1-k-i}$
for all $i$ and $h_{m+2-k-i}=h_{m+3-k-i}=\cdots=h_m=0.$\\
(5) $h_i \ge 0 $ for all $i$ (Nonnegativity)
(See Stanley \cite{ST1, ST3}.)\\
(6) If $G_1$ is a subgraph of $G_2$, then ${\mathcal P}(G_1) \subset {\mathcal P}(G_2)$\\
and $h_i(G_1) \le h_i(G_2).$ (Monotonicity)(See Stanley \cite{ST1, ST3}.)\\
(7) $volume({\mathcal P}(G))=\frac{P(1)}{m!}=$ leading coefficient of $WG(n).$ \\
(See Stanley \cite{ehrhart, ST4}.)\\
(8) All real roots $ \alpha $ of $WG(n)=0$ satisfy
 $-m \le \alpha < \lfloor \frac{m}{2} \rfloor.$(See Beck et al \cite{BDDPS}.)
\end{rmk}

\begin{rmk}
There seem to be some connections to {\bf semi magic cubes}. A
semi magic cube is an $k \times k \times \cdots \times k$ table
with nonnegative integral entries such that each directional sum
must be equal to $n$ (a {\bf magic sum}). (The reader can find
more details on semi magic squares and semi magic cubes in
\cite{BeckRobins}). For example, suppose we have a $2 \times 2$
magic square.  Then, the defining polytope for a $2 \times 2$
magic square with a magic sum $n$ is a face of the polytope
defining a weighted cycle graph $WC_4(n)$. In general suppose we
have a $2^d$ magic square.  Then the defining polytope for a $2^d$
magic square with a magic sum $n$ is a face of the polytope
defining a weighted $d-$dimensional Hypercubic Graph $HC_d$.
\end{rmk}

\begin{rmk}
We can get $\rho(G)$ from the information of the graph $G=(V,E)$
using the package {\sffamily LattE} by J.A.De Loera et. al.
(\cite{DHHHTY}), and the Elliott Maple package by G. Xin
(\cite{GX}), which improved the Omega package by G. Andrews et.
al. (\cite{APR} or references therein(in fact, their serial
articles)). We also can find the coordinates of all vertices of
the polytope $P(G)$ using {\sffamily cdd} and {\sffamily cdd+} by
K. Fukuda (\cite{FU}).

\end{rmk}
\section{Further Questions}\label{question}

\begin{enumerate}
\item It is obvious that, for a given number $m$ of vertices,
$$WK_m(n) \le WG(n) \le WD_m(n)$$ for any simple graph $G=(V,E)$ with $|V|=m,$
where $K_n$(resp. $D_m$) is a complete(resp. discrete ) graph of
order $m$. We can ask the same question for trees. That is, given
the number of vertices what kind of trees $T$ achieve the  maximal
or minimal values of $WT(n)$.

\item If a simple graph $G$ is {\bf connected} and
$\rho(G)=\frac{P(x)}{Q(x)}$, where $P(x)$ and $Q(x)$ are
polynomials of lowest degree, then is $P(x)$ always symmetric ?
That is, if $k=degree(P(x))$ then does $x^kP(\frac{1}{x})=P(x)$
hold?

\item For a given simple graph $G$, can we express the volume of
the polytope ${\mathcal P}(G)$ in terms of  the graph $G$?

\item If ${\bullet}$ is an operation between two simple graphs
$G_1,G_2$ and it is closed in the simple graphs, what is
$W(G_1{\bullet}G_2)(n)$ and $\rho(G_1{\bullet}G_2)$? Can we get
any  relation between $WG(n)$ (resp. $\rho(G)$) and $W{\bar G}(n)$
(resp. $\rho({\bar G})$)? (${\bar G}$ is a complement graph of
$G$).

\item Given two simple connected graphs $G_1=(V,E_1)$ and
$G_2=(W,E_2)$, choose a vertex $v \in V$,$w \in W$. Make $v$ and
$w$ adjacent by adding an edge between them so that two separate
graphs $G_1$ and $G_2$ becomes one connected graph $G$. Which
vertices in each side do we have to choose in order to maximize or
minimize $WG(n)$ ?

\item Generalization or Extension of the Weighted Graph by
inserting one or several slack variables in the edge between two
adjacent vertices(cf. Gear Graph), or to the Problems related to
(or symmetric) Magic Squares.\\
\xymatrix{i \ar@{-}[r] &{\bullet} \ar@{-}[r] &\cdots
\ar@{-}[r] &{\bullet} \ar@{-}[r] &j\\&&&&} \qquad $\Rightarrow\\
M=(f_{ij}(n,k)):(n+1)\times (n+1)$, where $$f_{ij}(n,k)=
\frac{\lfloor 2^{n-i-j}\rfloor}{2^{n-i-j}} C(n+k-i-j,k)$$ for $0
\le i,j \le n$. Note that $M=(f_{ij}(n,0))=B(n)$.

\item For what kinds of rational function $f(x)$ does there exist
a simple graph $G$ so that $\rho(G)=f(x)$? In other words, what is
the image $\rho({\mathcal G})$ of the map $\rho$ in ${\bf Z}[[x]]$
?

\item Compute $WG(n)$ and $\rho(G)$  for other shape of simple
graphs, like Wheel Graph, Cayley Graph, Complete k-partite Graph,
Web Graph or Grid Graph, Regular Graphs, other Platonic Graphs
etc...

\end{enumerate}


\end{document}